\documentclass[11pt]{amsart}
\usepackage{times,epsf,amssymb,amsmath,hyperref}

\begin{document}

\newtheorem{thm}{Theorem}[section]
\newtheorem{lem}[thm]{Lemma}
\newtheorem{cor}[thm]{Corollary}

\theoremstyle{definition}
\newtheorem{defn}{Definition}[section]

\theoremstyle{remark}
\newtheorem{rmk}{Remark}[section]

\def\square{\hfill${\vcenter{\vbox{\hrule height.4pt \hbox{\vrule
width.4pt height7pt \kern7pt \vrule width.4pt} \hrule height.4pt}}}$}

\def\T{\mathcal T}

\newenvironment{pf}{{\it Proof:}\quad}{\square \vskip 12pt}

\title{Number of Least Area Planes in Gromov Hyperbolic $3$-Spaces}
\author{Baris Coskunuzer}
\address{Koc University \\ Department of Mathematics \\ Sariyer, Istanbul 34450 Turkey}
\email{bcoskunuzer@ku.edu.tr}
\thanks{The author is partially supported by NSF Grant DMS-0603532 and TUBITAK Grant 107T642}

\maketitle


\newcommand{\Si}{S^2_{\infty}(\mathbf{H}^3)}
\newcommand{\PI}{\partial_{\infty}}
\newcommand{\SI}{S^2_{\infty}}
\newcommand{\SX}{S^2_{\infty}(X)}
\newcommand{\SM}{S^2_{\infty}(\widetilde{M})}
\newcommand{\BHH}{\mathbf{H}^3}
\newcommand{\CH}{\mathcal{C}(\Gamma)}
\newcommand{\BH}{\mathbf{H}}
\newcommand{\BR}{\mathbf{R}}
\newcommand{\BC}{\mathbf{C}}
\newcommand{\BZ}{\mathbf{Z}}

\begin{abstract}

We show that for a generic simple closed curve $\Gamma$ in the asymptotic boundary of a Gromov
hyperbolic $3$-space with cocompact metric $X$, there exist a unique least area plane $\Sigma$ in
$X$ such that $\PI \Sigma = \Gamma$. This result has interesting topological applications for
constructions of canonical $2$-dimensional objects in $3$-manifolds.

\end{abstract}

\section{Introduction}

Let $X$ be a Gromov hyperbolic $3$-space with cocompact metric, and $\SX$ be the sphere at infinity
of $X$. Let $\Gamma$ be a given simple closed curve in $\SX$. The asymptotic Plateau problem asks
the existence of a least area plane $\Sigma$ in $X$ with asymptotic boundary $\Gamma$, i.e. $\PI
\Sigma = \Gamma$. The author gave a positive answer to this question, and showed the existence of
such least area plane asymptotic to given curve in $\SX$ in \cite{Co1}.

For the analogous problem for $\BH^3$, Michael Anderson showed the existence of a least area plane
$\Sigma\subset\BH^3$ asymptotic to given simple closed curve $\Gamma\subset \Si$ in \cite{A1},
\cite{A2}. Later, by using topological techniques instead of geometric measure theory, Gabai got a
similar result for hyperbolic $3$-space with cocompact metric in \cite{Ga}.

On the other hand, on the number of the least area planes for a given simple closed curve in the
asymptotic sphere, there are a few results so far. In \cite{A1}, Anderson showed that if the given
asymptotic boundary $\Gamma\subset\Si$ bounds a convex domain in $\Si$, then there exists a unique
least area plane in $\BH^3$ . Then, Hardt and Lin generalized this result to the simple closed
curves bounding star shaped domains in $\Si$ in \cite{HL}. Recently, the author showed that a
generic simple closed curve in $\Si$ bounds a unique least area plane in $\BH^3$ in \cite{Co2},
\cite{Co4}.

For the number of least area planes in a Gromov hyperbolic spaces spanning a given curve in the
asymptotic boundary, there is no result so far. In this paper, we will generalize the results in
\cite{Co4} to Gromov hyperbolic $3$-spaces with cocompact metric. We will use similar techniques,
but there are more technicalities to overcome because of Gromov hyperbolic space \cite{Co3}. Our
main result is as follows:

\vspace{0.3cm}

\noindent \textbf{Theorem 3.3.} Let $X$ be a Gromov hyperbolic $3$-space with cocompact metric. Let
$A$ be the space of simple closed curves in $\SX$ and $A'\subset A$ be the subspace containing the
simple closed curves in $\SX$ bounding a unique least area plane in $X$. Then, $A'$ is generic in
$A$, i.e. $A-A'$ is a set of first category.

\vspace{0.3cm}

This is the first result in this direction for Gromov hyperbolic spaces. The techniques are very
general, and they can be applied to many similar settings. Furthermore, there are interesting
applications of this result. Namely, for a generic curve as above, the following conjecture is
true.

\vspace{0.3cm}

\noindent {\bf Disjoint Planes Conjecture:} Let $\Gamma_1, \Gamma_2$ be simple closed curves in
$\SX$. If $\Gamma_1$ and $\Gamma_2$ do not cross each other (i.e. They are the boundaries of
disjoint open regions in $\SX$ ), then any distinct least area planes $\Sigma_1, \Sigma_2$ in $X$
with asymptotic boundary $\Gamma_1, \Gamma_2$ are disjoint.

\vspace{0.3cm}

With this conjecture, it can be showed that for a genuine lamination $\Lambda$ in $M$, If
$\widetilde{\Lambda}$ has continuous extension property, then by replacing the leaves if
$\widetilde{\Lambda}$ with least area planes $\widetilde{\Delta}$, we get a genuine lamination
$\Delta$ in $M$ with minimal leaves. Hence, if the induced curves in $\SX$ are generic, with the
main result, it is possible to modify a given genuine $\Lambda$ to the one with minimal leaves
$\Delta$ in $M$. Of course, it is possible to do same for any essential $2$-dimensional object in a
Gromov hyperbolic $3$-manifold.

The organization of the paper is as follows: In the next section we will cover some basic results
which will be used in the following sections. In section 3, we will prove the main result of the
paper. Then in section 4, we will show some applications of the main result. Finally in section 5,
we will have some final remarks.

\section{Preliminaries}

A $3$-manifold $M$ is called {\em Gromov hyperbolic manifold} if its fundamental group $\pi_1(M)$
is a word hyperbolic (or Gromov hyperbolic) group \cite{Gr}. We call $X$ as {\em Gromov hyperbolic
$3$-space with cocompact metric} if $X$ is the universal cover of a Riemannian closed orientable
irreducible Gromov hyperbolic $3$-manifold $M$ where the metric on $X$ is induced by $M$. By
\cite{BM}, $X$ is homeomorphic to an open ball in $\BR^3$. Since $X$ is Gromov hyperbolic
$3$-space, it has a natural compactification $\overline{X}$ where $\overline{X}= X \cup \PI X$.
Here, $\PI X$ is the sphere at infinity $\SX$, and a point on $\SX$ corresponds to an equivalence
class of infinite rays in $X$ where two rays are equivalent if they are asymptotic.

We call a disk $D$ as a {\em least area disk} if $D$ has the smallest area among the disks with the
same boundary $\partial D$. A {\em plane} is a subset of $X$ with the topological type of a disk,
and complete with the induced path metric from $X$. We call a plane $P$ as a {\em least area plane}
if any subdisk in the plane $P$ is a least area disk.

A codimension-$1$ lamination $\sigma$ in $X$ is a foliation of a closed subset of $X$ with
$2$-manifolds (\textit{leaves}) such that $X$ is covered by charts of the form $I^2\times I$ where
a leave passes through a chart in a slice of the form $I^2\times \{p\}$ for $p\in I$. Here and
later, we abuse notation by letting $\sigma$ also denote the underlying space of its lamination.

The sequence $\{D_{i}\}$ of embedded disks in a Riemannian manifold $X$ \textit{converges} to the
lamination $\sigma$ if

i) For any convergent sequence $\{x_{n_i} \}$ in $X$ with  $x_{n_i} \in D_{n_i}$ where $n_i$ is a
strictly increasing sequence, $\lim x_{n_i} \in \sigma$.

ii) For any $x \in \sigma$, there exists a sequence $\{x_i\}$ with $x_i \in D_i$ and $\lim x_i = x$
such that there exist embeddings $f_{i}: D^2 \to D_i$ which converge in the $C^{\infty }$-topology
to a smooth embedding $f:D^2 \to L_{x}$, where $x_i \in f_{i}(Int(D^2))$, and $L_{x}$ is the leaf
of $\sigma $ through $x$, and $x\in f(Int(D^2))$.

We call such a lamination $\sigma$ as \textit{$D^2$-limit lamination}.

\begin{thm}\cite{Co1}
Let $\Gamma $ be a simple closed curve in $\SX$ where $X$ is a Gromov hyperbolic $3$-space with
cocompact metric.  Then there exists a $D^2$-limit lamination $\sigma \subset X$ by least area
planes spanning $\Gamma$.
\end{thm}

\begin{rmk}
Gabai proved this theorem for hyperbolic $3$-spaces with cocompact metric in \cite{Ga}.
\end{rmk}

The following lemma is often called Meeks-Yau exchange roundoff trick.

\begin{lem}\cite{MY} If $D_1$ and $D_2$ are least area disks in a manifold $M$ such that $D_1\cap\partial D_2 =
\emptyset$ and $D_2\cap\partial D_1 = \emptyset$, then $D_1 \cap D_2=\emptyset$.
\end{lem}

\section{Generic Uniqueness of Least Area Planes}

In this section, we will prove that the space of simple closed curves in $\SX$ bounding a unique
least area plane in $X$ is generic in the space of simple closed curves in $\SX$.

The short outline of the technique is the following: Let $\Gamma_0$ be a simple closed curve in
$\SX$. First, we will show that either there exists a unique least area plane $\Sigma_0$ in $X$
with $\PI\Sigma_0=\Gamma_0$, or there exist two {\em disjoint} least area planes $\Sigma_0^+ ,
\Sigma_0^-$ in $X$ with $\PI\Sigma_0^\pm=\Gamma_0$.

Now, take a small neighborhood $N(\Gamma_0)\subset \SX$ which is an annulus. Then foliate
$N(\Gamma_0)$ by simple closed curves $\{\Gamma_t\}$ where $t\in(-\epsilon, \epsilon)$, i.e.
$N(\Gamma_0) \simeq \Gamma\times (-\epsilon, \epsilon)$. By the above fact, for any $\Gamma_t$
either there exist a unique least area plane $\Sigma_t$, or there are two least area planes
$\Sigma_t^\pm$ disjoint from each other. Also, since these are least area planes, if they have
disjoint asymptotic boundary, then they are disjoint by Meeks-Yau exchange roundoff trick. This
means, if $t_1<t_2$, then $\Sigma_{t_1}$ is disjoint and \textit{below} from $\Sigma_{t_2}$ in $X$.
Consider this collection of least area planes. Note that for curves $\Gamma_t$ bounding more than
one least area plane, we have a canonical region $N_t$ in $X$ between the disjoint least area
planes $\Sigma_t^\pm$.

Now, $N(\Gamma)$ separates $\SX$ into two parts, and take a proper line $\beta\subset X$ which is
asymptotic to two points belongs to these two different parts. This line is transverse to the
collection of these least area planes asymptotic to the curves in $\{\Gamma_t\}$. Also, a finite
segment of this line intersects the entire collection. Let the length of this finite segment be
$C$.

Now, the idea is to consider the {\em thickness} of the neighborhoods $N_t$ assigned to the
asymptotic curves $\{\Gamma_t\}$. Let $s_t$ be the length of the segment $I_t$ of $\beta$ between
$\Sigma_t^+$ and $\Sigma_t^-$, which is the {\em width} of $N_t$ assigned to $\Gamma_t$. Then, the
curves $\Gamma_t$ bounding more than one least area planes have positive width, and contributes to
total thickness of the collection, and the curves bounding unique least area plane has $0$ width
and do not contribute to the total thickness. Since $\sum_{t\in(-\epsilon, \epsilon)} s_t < C$, the
total thickness is finite. This implies for only countably many $t\in(-\epsilon, \epsilon)$,
$s_t>0$, i.e. $\Gamma_t$ bounds more than one least area plane. For the remaining uncountably many
$t\in(-\epsilon, \epsilon)$, $s_t=0$, and there exist a unique least area plane for those $t$. This
proves the space of Jordan curves of uniqueness is dense in the space of Jordan curves in $\SX$.
Then, by using similar arguments, we will show this space is not only dense, but also generic.

First, we will show that if two least area planes have disjoint asymptotic boundaries, then they
are disjoint.

\begin{lem}
Let $\Gamma_1$ and $\Gamma_2$ be two disjoint simple closed curves in $\SX$. If $\Sigma_1$ and
$\Sigma_2$ are least area planes in $X$ with $\PI \Sigma_i = \Gamma_i$, then $\Sigma_1$ and
$\Sigma_2$ are disjoint, too.
\end{lem}

\begin{pf}Assume that $\Sigma_1\cap\Sigma_2\neq\emptyset$. So, the intersection between
$\Sigma_1$ and $\Sigma_2$ contains a simple closed curve or infinite line.

Let's assume the intersection contains a simple closed curve $\gamma$. Since $\Sigma_1$ and
$\Sigma_2$ are both minimal, the intersection must be transverse on a subarc of $\gamma$ by maximum
principle. Now, $\gamma$ bounds two least area disks $D_1$ and $D_2$ in $X$, with
$D_i\subset\Sigma_i$. Now, take a larger subdisk $E_1$ of $\Sigma_1$ containing $D_1$, i.e.
$D_1\subset E_1 \subset \Sigma_1$. By definition, $E_1$ is also an least area disk. Now, modify
$E_1$ by swapping the disks $D_1$ and $D_2$. Then, we get a new disk $E_1 '= \{E_1 - D_1\} \cup
D_2$. Now, $E_1$ and $E_1 '$ have same area, but $E_1 '$ have a folding curve along $\gamma$. By
smoothing out this curve as in \cite{MY} (Meeks-Yau exchange roundoff trick), we get a disk with
smaller area, which contradicts to $E_1$ being least area.

If the intersection contains an infinite line $l$, since asymptotic boundaries $\Gamma_1$ and
$\Gamma_2$ are disjoint, $l$ cannot limit on the $\SI(X)$. So, $l$ must be in the bounded subset of
$X$. If $\Sigma_1$ or $\Sigma$ is properly embedded, they cannot contain an infinite line in the
intersection. This is because, in that case, the limit point must be in the intersection, and this
contradicts to the characterization of intersection of minimal surfaces, \cite{HS}. So, $\Sigma_1$
and $\Sigma_2$ are nonproperly embedded least area planes. Same idea shows that $l$ is a proper
line in $\Sigma_i$ with the induced path metric.

Let $\{D_i\}$ be a sequence of subdisks in $\Sigma_1$ such that $D_1\subset D_2 \subset ... \subset
D_n \subset ... $ with $\partial D_i \rightarrow \PI \Sigma_1$. By \cite{Co1}, $\{D_i\}$ induces a
$D^2$-limit lamination $\widehat{\Sigma}_1$ with $\PI \widehat{\Sigma}_1 = \PI \Sigma_1$. This
lamination is not empty since $\Sigma_1$ is a leaf in the lamination by construction. Similarly,
define $\widehat{\Sigma}_2$.

Now, since $\widehat{\Sigma}_i$ is closed subset of $X$, the limit set of $l$ must be in
$\widehat{\Sigma}_1 \cap \widehat{\Sigma}_2$. We claim that the limit set of $l$ produce a simple
closed curve in the intersection. Let $p$ be a limit point of $l$, and it is contained in the leaf
$\Sigma_1'$ and $\Sigma_2'$ where $\Sigma_i'$ is a least area plane in $\widehat{\Sigma}_i$. By the
definition of $D^2$-limit lamination, there is a small disk $E$ in $\Sigma_1'$ containing $p$ such
that $E_i \rightarrow E$ in $C^\infty$ topology where each $E_i$ is a small disk in
$D_i\subset\Sigma_1$. Since $p$ is in the limit set of $l$, by the convergence of $\{E_i\}$, $p$
must be an interior point in the limit set of $l$. This implies $\Sigma_1'\cap\Sigma_2'$ contains a
simple closed curve, or a proper line in the induced metric. However, containing an infinite line
in the limit implies the existence of monogons in $\Sigma_1$, and $\Sigma_2$ which contradicts
being least area, \cite{HS}. So, this implies the existence of a simple closed curve in the
intersection of least area planes $\Sigma_1'$ and $\Sigma_2'$. Again, this is a contradiction by
the second paragraph of this proof.

\end{pf}

The following lemma is very essential for our technique. Mainly, the lemma says that for any given
simple closed curve $\Gamma$ in $\SX$, either there exists a unique least area plane $\Sigma$ in
$X$ asymptotic to $\Gamma$, or there exist two least area planes $\Sigma^\pm$ in $X$ which are
asymptotic to $\Gamma$ and disjoint from each other.

\begin{lem}
Let $X$ be a Gromov hyperbolic $3$-space with cocompact metric, and let $\Gamma$ be a simple closed
curve in $\SX$. Then either there exists a unique least area plane $\Sigma$ in $X$ with
$\PI\Sigma=\Gamma$, or there are two canonical disjoint extremal least area planes $\Sigma^+$ and
$\Sigma^-$ in $X$ with $\PI \Sigma^\pm = \Gamma$. Moreover, any least area plane $\Sigma '$ with
$\PI\Sigma '= \Gamma$ is disjoint from $\Sigma^\pm$, and it is captured in the region bounded by
$\Sigma^+$ and $\Sigma^-$ in $X$.
\end{lem}

\begin{pf} Let $\Gamma$ be a simple closed curve in $\SX$. $\Gamma$ separates $\SX$ into two parts, say
$\Omega^+$ and $\Omega^-$ which are open disks. Define sequences of pairwise disjoint simple closed
curves $\{\Gamma_i^+\}$  such that $\Gamma_i^+ = \partial E_i$ where $E_i$ is a closed disk in
$\Omega^+$ for any $i$, and $E_i\subset Int(E_j)$ for any $i<j$. Moreover, $\Omega^+ =\bigcup_i
E_i$ and $\Gamma_i^+ \rightarrow \Gamma$.

By Theorem 2.1, for any $\Gamma_i^+ \subset \SX$, there exist a least area plane $\Sigma_i^+$ with
$\PI \Sigma_i^+ = \Gamma_i^+$. Note that these least area planes are separating in $X$ by their
construction in Theorem 2.1. Also, let $p^+ \in Int(E_1) \subset \Omega^+$ and $p^- \in \Omega^-$
be two points belonging different components of $\Si-\Gamma$. Let $\beta$ be a proper line in $X$
asymptotic to $p^+$ and $p^-$. Since each $\overline{\Sigma_i^+}$ is separating in $\overline{X}$,
$\beta$ intersects $\Sigma_i^+$ for any $i$. Let $x_i$ be a point in $\Sigma_i^+ \cap \beta$ for
any $i$.

Now, we define the sequence of least area disks. Let $D_i$ be a least area disk in $\Sigma_i^+$
containing $x_i$ such that $i< d_X(x_i,\partial D_i) < i+1$ where $d_X$ is the extrinsic distance
in $X$. We claim that the sequence of least area disks $\{D_i\}$ in $X$ converge to a nonempty
lamination $\sigma$.  As $\alpha_i =\partial D_i$ converges to $\Gamma$, by Theorem 2.1, $D_i$
converges to a lamination $\sigma$ (possibly empty) by least area planes with $\PI \sigma =\Gamma$.
Now, we show that $\sigma$ is a nonempty lamination by least area planes. Since $p^+ \in \Omega^+$
and $p^-\in \Omega^-$ where $\Omega^\pm$ are open disks in $\SX$, we can find a sufficiently small
$\epsilon$ such that $B_\epsilon(p^+) \subset Int(E_1)$ and $B_\epsilon(p^-) \subset \Omega^-$. Let
$\gamma^+ = \partial B_\epsilon(p^+)$ and $\gamma^- = \partial B_\epsilon(p^-)$. Then, $\gamma^\pm$
are two simple closed curves disjoint from $\Gamma$. By Theorem 2.1, there are least area planes
$P^+$ and $P^-$ in $X$ with $\PI P^\pm = \gamma^\pm$. By Lemma 3.1, $P^\pm$ are disjoint from
$\Sigma_i^+$ for each $i$. So, the finite segment $\hat{\beta}$ of $\beta$ between $P^+$ and $P^-$
contains all intersection points $\Sigma_i^+ \cap \beta$ for any $i$. Hence,
$\{x_i\}\subset\hat{\beta}$. Since $\hat{\beta}$ is finite segment, it is compact, and the sequence
$\{x_i\}$ has a limit point $x$. This shows that the sequence $\{D_i\}$ has nonempty limit. So,
$\sigma$ is a nonempty lamination by least area planes with $\PI \sigma = \Gamma$ as claimed.

Now, we want to show that $\sigma$ consists of only one least area plane $\mathcal{P}$, i.e.
$\sigma = \mathcal{P}$. Assume that $\sigma$ contains more than one least area plane. Recall that
$\sigma$ is a collection of disjoint least area planes asymptotic to $\Gamma$, and $\sigma$ is a
closed subset of $X$. Note that each plane $\mathcal{P}$ in $\sigma$ is separating in $X$ by
construction by Theorem 4.3 Step $4$ of \cite{Co1} (or Proposition 3.9 Step $3$ in \cite{Ga}).
Consider the components of $X- \sigma$. Let $X^+$ be the component of $X - \sigma$ with $\PI X^+ =
\Omega^+$. Note that as $\sigma$ is closed, each component of $X-\sigma$ is open, and so is $X^+$.
Since for each $i$, $D_i \subset \Sigma_i^+$ and $\PI \Sigma_i^+ = \Gamma_i^+ \subset \Omega^+$, by
Lemma 3.2, $D_i \subset X^+$. Since each plane $\mathcal{P}$ in $\sigma$ is separating, there is a
plane $\mathcal{P}_1$ such that $\mathcal{P}_1 = \partial X^+$. Let $\mathcal{P}_2$ be another
least area plane in $\sigma$. We claim that there cannot be such a plane $\mathcal{P}_2$ because of
the special properties of the sequence of least area disks $\{D_i\}$. Let $Y^+$ be the component of
$X-\mathcal{P}_2$ with $\PI Y^+ = \Omega^+$. Clearly, $X^+ \subset Y^+$. Moreover, since either
there exist an open complementary region (a component of $X-\sigma$), or there is an open region
foliated by least area planes in $\sigma$ between $\mathcal{P}_1$ and $\mathcal{P}_2$, hence $X^+
\subsetneq Int(Y^+)$. This means $\mathcal{P}_1$ forms a barrier between the sequence of least area
disks $\{D_i\}$ and $\mathcal{P}_2$. In other words, the sequence $\{D_i\}$ cannot reach
$\mathcal{P}_2$, so $\mathcal{P}_2$ cannot be in the limit. More precisely, if $q \in
\mathcal{P}_2$, since $D_i\subset X^+$ for any $i$, and $X^+ \subsetneq Int(Y^+)$, there is no
sequence $\{q_i\}$ with $q_i \in D_i$ such that $q_i \to q$. So, $\mathcal{P}_2$ cannot be in the
limit of $\{D_i\}$. Hence, this shows $\sigma$ consists of only one plane $\mathcal{P}_1$, i.e.
$\sigma = \mathcal{P}_1$. Call this least area plane $\mathcal{P}_1$ as $\Sigma^+$. Similarly by
defining a similar sequence of simple closed curves converging to $\Gamma$ in $\Omega^-$, we can
define the least area plane $\Sigma^+$.

Now, let $\Sigma '$ be any other least area plane in $X$ with $\PI\Sigma '=\Gamma$. If $\Sigma
'\cap\Sigma^+ \neq \emptyset$, then some part of $\Sigma '$ must be {\em above} $\Sigma^+$. Since
$\Sigma^+ =\lim D_i$ where $D_i\subset \Sigma_i^+$, for sufficiently large $i$, $\Sigma '
\cap\Sigma_i^+ \neq\emptyset$. However, $\PI\Sigma_i^+ = \Gamma_i^+$ is disjoint from
$\Gamma=\PI\Sigma '$. Then, by Lemma 3.1, $\Sigma '$ must be disjoint from $\Sigma_i^+$. This is a
contradiction.

Similarly, this is true for $\Sigma^-$, too. Moreover, let $N \subset X$ be the region between
$\Sigma^+$ and $\Sigma^-$, i.e. $\partial N= \Sigma^+ \cup \Sigma^-$. Then by construction,  $N$ is
also a canonical region for $\Gamma$, and for any least area plane $\Sigma '$ with $\PI \Sigma
'=\Gamma$, $\Sigma '$ is contained in the region $N$, i.e. $\Sigma '\subset N$. This shows that if
$\Sigma^+ = \Sigma^-$, there exist a unique least area plane asymptotic to $\Gamma$. If $\Sigma^+
\neq \Sigma^-$, then they must be disjoint.
\end{pf}

Now, we are going to prove the main theorem of the paper. This theorem says that for a generic
simple closed curve in the asymptotic sphere of a Gromov hyperbolic $3$-space, there exists a
unique least area plane.

\begin{thm} Let $X$ be a Gromov hyperbolic $3$-space with cocompact metric, and let
$A$ be the space of simple closed curves in $\SX$. Let $A'\subset A$ be the subspace containing the
simple closed curves in $\SX$ bounding a unique least area plane in $X$. Then, $A'$ is generic in
$A$, i.e. $A-A'$ is a set of first category.
\end{thm}

\begin{pf} We will prove this theorem in 2 steps.

\vspace{0.3cm}

\textbf{Claim 1:} $A'$ is dense in $A$ as a subspace of $C^0(S^1,\SX)$ with the supremum metric.

\vspace{0.3cm}

\begin{pf} $A$ is the space of Jordan curves in $\SX$.\\
Then, $A= \{\alpha\in C^0(S^1,S^2)\ | \ \alpha(S^1)$ is an embedding $\}$.

Now, let $\Gamma_0\in A$ be a simple closed curve in $\SX$. Since $\Gamma_0$ is simple, there exist
a small neighborhood $N(\Gamma_0)$ of $\Gamma_0$ which is an annulus in $\SX$. Let
$\Gamma:(-\epsilon,\epsilon)\rightarrow A$ be a small path in $A$ through $\Gamma_0$ such that
$\Gamma(t)=\Gamma_t$ and $\{\Gamma_t\}$ foliates $N(\Gamma)$ with simple closed curves $\Gamma_t$.
In other words, $\{\Gamma_t\}$ are pairwise disjoint simple closed curves, and
$N(\Gamma_0)=\bigcup_{t\in (-\epsilon,\epsilon)} \Gamma_t$.

Now, $N(\Gamma_0)$ separates $\SX$ into two parts, say $D^+$ and $D^-$, i.e. $\SX=N(\Gamma_0)\cup
D^+\cup D^-$. Let $p^+$ be a point in $D^+$ and let $p^-$ be a point in $D^-$ such that for a small
$\delta$, $B_\delta(p^\pm)$ are in the interior of $D^\pm$. Let $\beta$ be a proper line in $X$
asymptotic to $p^+$ and $p^-$.

By Lemma 3.2, for any $\Gamma_t$ either there exist a unique least area plane $\Sigma_t$ in $X$, or
there is a canonical region $N_t$ in $X$ between the canonical least area planes $\Sigma_t^+$ and
$\Sigma_t^-$. With abuse of notation, if $\Gamma_t$ bounds a unique least area plane $\Sigma_t$ in
$X$, define $N_t=\Sigma_t$ as a degenerate canonical neighborhood for $\Gamma_t$. Then, let
$\widehat{N}= \{N_t\}$ be the collection of these degenerate and nondegenerate canonical
neighborhoods for $t\in(-\epsilon,\epsilon)$. Clearly, degenerate neighborhood $N_t$ means
$\Gamma_t$ bounds unique least area plane, and nondegenerate neighborhood $N_s$ means that
$\Gamma_s$ bounds more than one least area plane. Note that by Lemma 3.1, all canonical
neighborhoods in the collection are pairwise disjoint. On the other hand, by construction the
proper line $\beta$ intersects all the canonical neighborhoods in the collection $\widehat{N}$.

We claim that the part of $\beta$ which intersects $\widehat{N}$ is a finite line segment. Let
$P^+$ be the least area plane asymptotic to the round circle $\partial B_\delta(p^+)$ in $D^+$.
Similarly, define $P^-$. By Lemma 3.1, $P^\pm$ are disjoint from the collection of canonical
regions $\widehat{N}$. Let $\beta\cap P^\pm=\{q^\pm\}$. Then the part of $\beta$ which intersects
$\widehat{N}$ is the line segment $l\subset \beta$ with endpoints $q^+$ and $q^-$. Let $C$ be the
length of this line segment $l$.

Now, for each $t\in(-\epsilon,\epsilon)$, we will assign a real number $s_t\geq 0$. If there exists
a unique least area plane $\Sigma_t$ in $X$ for $\Gamma_t$ ($N_t$ is degenerate), then let $s_t$ be
$0$. If not, let $I_t = \beta\cap N_t$, and $s_t$ be the length of $I_t$. Clearly if $\Gamma_t$
bounds more than one least area plane ($N_t$ is nondegenerate), then $s_t > 0$. Also, it is clear
that for any $t$, $I_t\subset l$ and $I_t\cap I_s=\emptyset$ for any $t\neq s$. Then,
$\sum_{t\in(-\epsilon,\epsilon)} s_t < C$ where $C$ is the length of $l$. This means for only
countably many $t\in(-\epsilon,\epsilon)$, $s_t > 0$. So, there are only countably many
nondegenerate $N_t$ for $t\in(-\epsilon,\epsilon)$. Hence, for all other $t$, $N_t$ is degenerate.
This means there exist uncountably many $t\in(-\epsilon,\epsilon)$, where $\Gamma_t$ bounds a
unique least area plane. Since $\Gamma_0$ is arbitrary, this proves $A '$ is dense in $A$.
\end{pf}

\textbf{Claim 2:} $A'$ is generic in $A$, i.e. $A-A'$ is a set of first category.

\vspace{0.3cm}

\begin{pf} We will prove that $A '$ is countable intersection of open dense subsets of a
complete metric space. Then the result will follow by Baire category theorem.

Since the space of continuous maps from circle to sphere $C^0(S^1,S^2)$ is complete with supremum
metric, then the closure of $A$ in $C^0(S^1,S^2)$, $\bar{A}\subset C^0(S^1,S^2)$, is also complete.

Now, we will define a sequence of open dense subsets $U^i\subset A$ such that their intersection
will give us $A '$. Let $\Gamma\in A$ be a simple closed curve in $\SX$, as in the Claim 1. Let
$N(\Gamma)\subset \SX$ be a neighborhood of $\Gamma$ in $\SX$, which is an open annulus. Then,
define an open neighborhood $U_\Gamma$ of $\Gamma$ in $A$, such that $U_\Gamma = \{\alpha \in A \ |
\ \alpha(S^1)\subset N(\Gamma), \ \alpha \mbox{ is homotopic to } \Gamma\}$. Clearly, $A=
\bigcup_{\Gamma\in A} U_\Gamma$.  Now, define a proper line $\beta_\Gamma$ as in Claim 1, which
intersects all the least area planes asymptotic to curves in $U_\Gamma$.

Now, for any $\alpha \in U_\Gamma$, by Lemma 3.2, there exist a canonical region $N_\alpha$ in $X$
(which can be degenerate if $\alpha$ bounds a unique least area plane). Let $I_{\alpha,\Gamma} =
N_\alpha \cap \beta_\Gamma$. Then let $s_{\alpha,\Gamma}$ be the length of $I_{\alpha,\Gamma}$
($s_{\alpha,\Gamma}$ is $0$ if $N_\alpha$ degenerate). Hence, for every element $\alpha$ in
$U_\Gamma$, we assign a real number $s_{\alpha,\Gamma} \geq 0$.

Now, we define the sequence of open dense subsets in $U_\Gamma$. Let $U^i_\Gamma = \{\alpha\in
U_\Gamma \ | \ s_{\alpha,\Gamma} < 1/i \  \}$. We claim that $U^i_\Gamma$ is an open subset of
$U_\Gamma$ and $A$. Let $\alpha\in U^i_\Gamma$, and let $s_{\alpha,\Gamma} = \lambda < 1/i$. So,
the interval $I_{\alpha,\Gamma}\subset \beta_\Gamma$ has length $\lambda$. let $I ' \subset
\beta_\Gamma$ be an interval containing $I_{\alpha,\Gamma}$ in its interior, and has length less
than $1/i$. By the proof of Claim 1, we can find two simple closed curves $\alpha^+, \alpha^- \in
U_\Gamma$ with the following properties.

\begin{itemize}

\item $\alpha^\pm$ are disjoint from $\alpha$,

\item $\alpha^\pm$ are lying in opposite sides of $\alpha$ in $\SX$,

\item $\alpha^\pm$ bounds unique least area planes $\Sigma_{\alpha^\pm}$,

\item $\Sigma_{\alpha^\pm} \cap \beta_\Gamma \subset I '$.

\end{itemize}

The existence of such curves is clear from the proof Claim 1, as if one takes any foliation
$\{\alpha_t\}$ of a small neighborhood of $\alpha$ in $\SX$, there are uncountably many curves in
the family bounding a unique least area plane, and one can choose sufficiently close pair of curves
to $\alpha$, to ensure the conditions above.

After finding $\alpha^\pm$, consider the open annulus $F_\alpha$ in $\SX$ bounded by $\alpha^+$ and
$\alpha^-$. Let $V_\alpha = \{ \gamma\in U_\Gamma \ | \ \gamma(S^1)\subset F_\alpha , \ \gamma
\mbox{ is homotopic to } \alpha \}$. Clearly, $V_\alpha$ is an open subset of $U_\Gamma$. If we can
show $V_\alpha\subset U^i_\Gamma$, then this proves $U^i_\Gamma$ is open for any $i$ and any
$\Gamma\in A$.

Let $\gamma\in V_\alpha$ be any curve, and $N_\gamma$ be its canonical neighborhood given by Lemma
3.2. Since $\gamma(S^1)\subset F_\alpha$, $\alpha^+$ and $\alpha^-$ lie in opposite sides of
$\gamma$ in $\SX$. This means $\Sigma_{\alpha^+}$  and $\Sigma_{\alpha^-}$ lie in opposite sides of
$N_\gamma$. By choice of $\alpha^\pm$, this implies $N_\gamma \cap \beta_\Gamma= I_{\gamma,\Gamma}
\subset I '$. So, the length $s_{\gamma,\Gamma}$ is less than $1/i$. This implies $\gamma\in
U^i_\Gamma$, and so $V_\alpha\subset U^i_\Gamma$. Hence, $U^i_\Gamma$ is open in $U_\Gamma$ and
$A$.

Now, we can define the sequence of open dense subsets. let $U^i = \bigcup_{\Gamma\in A} U^i_\Gamma$
be an open subset of $A$. Since, the elements in $A '$ represent the curves bounding a unique area
minimizing plane, for any $\alpha\in A '$, and for any $\Gamma\in A$, $s_{\alpha,\Gamma} = 0$. This
means $A'\subset U^i$ for any $i$. By Claim 1, $U^i$ is open dense in $A$ for any $i>0$.

As we mention at the beginning of the proof, since the space of continuous maps from circle to
sphere $C^0(S^1,S^2)$ is complete with supremum metric, then the closure $\bar{A}$ of $A$ in
$C^0(S^1,S^2)$ is also complete metric space. Since $A'$ is dense in $A$, it is also dense in
$\bar{A}$. As $A$ is open in $C^0(S^1,S^2)$, this implies $U^i$ is a sequence of open dense subsets
of $\bar{A}$. On the other hand, since $s_{\alpha,\Gamma} = 0$ for any $\alpha\in A '$, and for any
$\Gamma\in A$, $A ' = \bigcap_{i>0} U^i$. Then, $A-A'$ is a set of first category, by Baire
Category Theorem. Hence, $A'$ is generic in $A$.
\end{pf}
\end{pf}

Hence, this proves that for a generic simple closed curve in the asymptotic sphere of a Gromov
hyperbolic $3$-space with cocompact metric, there exist a unique least area plane spanning the
curve. In the next section, we will give some applications of this result.

\section{Applications}

In this section, we will show that, if a simple closed curve $\Gamma$ in $\SX$, bounds a unique
least area plane $\Sigma_\Gamma$ in $X$, where $X$ is a Gromov hyperbolic $3$-space with cocompact
metric, then for any other simple closed curve $\alpha$ not crossing $\Gamma$ in $SX$, then any
least area planes $\Sigma_\alpha$ is disjoint from $\Sigma_\Gamma$. In other words, the following
Disjoint Planes Conjecture is true if one of the curves bounds a unique least area plane.

\vspace{0.3cm}

\noindent {\bf Disjoint Planes Conjecture:} Let $\Gamma_1, \Gamma_2$ be simple closed curves in
$\SX$, where $X$ is a Gromov hyperbolic $3$-space with cocompact metric. If $\Gamma_1$ and
$\Gamma_2$ do not cross each other (i.e. They are the boundaries of disjoint open regions in $\SX$
), then any distinct least area planes $\Sigma_1, \Sigma_2$ in $X$ with asymptotic boundary
$\Gamma_1, \Gamma_2$ are disjoint.

\vspace{0.3cm}

\begin{rmk}Even though this conjecture is interesting in its own right, it has powerful topological applications.
The most important application is constructing the least area representative of a $2$-dimensional
object in a $3$-manifold, like incompressible surfaces, and genuine laminations. For example, if
the $3$-manifold $M$ has a genuine lamination $\Lambda$, then by \cite{GK}, $M$ is automatically
Gromov hyperbolic. Then, up to the continuous extension property \cite{Fe}, \cite{Ca},
$\widetilde{\Lambda}$ will induce a $\pi_1$-invariant family of circles in $\SM$. Then by spanning
one of the curves in the family with a least area plane $\Sigma$, and considering images of
$\Sigma$ deck transformations, we get a $\pi_1$-invariant family of least area planes in $\SM$. By
the above conjecture, this family is pairwise disjoint, hence it can be showed that it projects
down to a embedded genuine lamination with least area leaves in $M$ by using the techniques in
\cite{Co1}. Similar construction works for incompressible surfaces, too, which implies a similar
result of \cite{HS}. In other words, with this conjecture, one can construct least area
representative of any essential $2$-dimensional object in a Gromov hyperbolic $3$-manifold.
\end{rmk}

Now, we will prove that if one of the curves in the conjecture bounds a unique least area plane,
then the conjecture is true.

\begin{thm} Let $\Gamma_1, \Gamma_2$ be as in the Disjoint Planes Conjecture. Let $\Gamma_1$ bounds
a unique least area plane in $X$. Then, the Disjoint Planes Conjecture is true for $\Gamma_1$ and
$\Gamma_2$.
\end{thm}

\begin{pf} First, we will show that, for any curves $\Gamma_1,\Gamma_2$ as in the statement of the theorem,
if $\Sigma_1^\pm$ and $\Sigma_2^\pm$ are the extremal least area planes with
$\Sigma_i^\pm=\Gamma_i$ as in the Lemma 3.2, then $\Sigma_1^+\cap\Sigma_2^+ = \emptyset$, and
likewise $\Sigma_1^-\cap\Sigma_2^- = \emptyset$.

Let $\SX-\Gamma_i=\Omega_i^+\cup\Omega_i^-$, and with loss of generality, further assume that
$\Omega_2^+\subset \Omega_1^+$ and $\Omega_1^-\subset \Omega_2^-$ (In a sense, $\Gamma_1$ {\em
stays below} of $\Gamma_2$ in $\SX$.) By the construction in Lemma 3.2, we know that
$\Sigma_1^+=\lim D_{1i}$ for some sequence of least area disks such that $D_{1i}\subset
\Sigma_{1i}$ where $\Sigma_{1i}$ is a least area plane with $\PI \Sigma_{1i} = \Gamma_{1i}$ and
$\Gamma_{1i}\subset\Omega_1^+$. Similarly same is true for $\Sigma_2^+$.

Assume that $\Sigma_1^+\cap\Sigma_2^+\neq\emptyset$. Then, by \cite{MY}, the intersection cannot
contain a simple closed curve. Hence, the intersection must contain an infinite line $l$ with $\PI
l \subset \Gamma_1\cap\Gamma_2$. This implies that some part of $\Sigma_1^+$ must be {\em above}
$\Sigma_2^+$. Recall that $\Sigma_2^+ =\lim D_{2i}$ where $D_{2i}\subset \Sigma_{2i}$. Moreover,
$\PI \Sigma_{2i} = \Gamma_{2i}$ and $\Gamma_{2i}\subset\Omega_2^+ \subset \Omega_1^+$. This means
for any $i$, $\Gamma_1\cap\Gamma_{2i} = \emptyset$. By Lemma 3.1, this implies
$\Sigma_1^+\cap\Sigma_{2i} = \emptyset$. However, as some part of $\Sigma_1^+$ is {\em above}
$\Sigma_2^+$, and $\Sigma_2^+ =\lim D_{2i}$ where $D_{2i}\subset \Sigma_{2i}$, this is a
contradiction. Hence $\Sigma_1^+\cap\Sigma_2^+ = \emptyset$. Similarly, $\Sigma_1^-\cap\Sigma_2^- =
\emptyset$ is also true for any such $\Gamma_1, \Gamma_2\subset \SX$.

Now, if $\Gamma_1$ is a simple closed curve in $\SX$ bounding a unique least area plane $\Sigma_1$
in $X$, then for $\Gamma_1$, the extremal least area planes $\Sigma_1^+$ and $\Sigma_2^+$ must be
same, i.e. $\Sigma_1^+=\Sigma_1^-=\Sigma_1$. Then, if $\Gamma_2$ is any other simple closed curve
in $SX$ which is not crossing $\Gamma_1$ in $SX$, then by previous paragraph,
$\Sigma_1^+\cap\Sigma_2^+ = \emptyset$ and $\Sigma_1^-\cap\Sigma_2^- = \emptyset$. If $\Sigma_2'$
is any least area plane with $\PI \Sigma_2' = \Gamma_2$, then by Lemma 3.2, $\Sigma_2'$ must be
{\em between} $\Sigma_2^+$ and $\Sigma_2^-$. Hence, as $\Sigma_1^+=\Sigma_1^-=\Sigma_1$ and
$\Sigma_1^\pm \cap \Sigma_2^\pm = \emptyset$, $\Sigma_1 \cap \Sigma_2' = \emptyset$. The proof
follows.
\end{pf}

By combining the above result with Theorem 3.3, we get the following corollary.

\begin{cor} The Disjoint Planes Conjecture is generically true.
\end{cor}

Now, by using the techniques in \cite{Co1}, as described in Remark 4.1, for a genuine lamination
$\Lambda$ in a $3$-manifold $M$, with the corollary above, it is generically possible to modify it
without changing the metric of $M$ so that the leaves of the new lamination $\Lambda '$ are
minimal. Similar, construction works for incompressible surfaces in Gromov hyperbolic
$3$-manifolds, to. Basically, this corollary implies that it is generically possible to construct
least area representative of any essential $2$-dimensional object in a Gromov hyperbolic
$3$-manifold.

\section{Final Remarks}

In this paper, we showed that for a generic simple closed curve in the asymptotic sphere of a
Gromov hyperbolic $3$-space with cocompact metric, there is a unique least area plane spanning the
curve in the space. This is the first result in this setting, while for the analogous question in
$\BH^3$, there are many significant results like \cite{A2}, \cite{HL}, \cite{Co2}, \cite{Co4}.

The next question is whether there exist simple closed curves in $\SX$ which bound more than one
least area planes in $X$, where $X$ is a Gromov hyperbolic $3$-space with cocompact metric. In the
special case, when $X=\BH^3$, we showed the existence of such curves in $\Si$ bounding more than
one least area plane in $\BH^3$ in \cite{Co5}. By using the techniques in \cite{Co5}, one can
easily show that if there exist a complete minimal surface $S$ with some genus in $X$ asymptotic to
a simple closed curve $\gamma_0$ in $\SX$, then there is a curve $\gamma_{t_0}$ in $\SX$ bounding
more than one least area plane in $X$ where $\{\gamma_t\}$ is a foliation of $\SX$ with simple
closed curves, i.e. $\SX=\bigcup_{t=-1}^1 \gamma_t$ ($\gamma_{-1}, \gamma_1$ are singular leaves).
The idea is by assuming the uniqueness for all curves, one can foliate $X$ with $\{\Sigma_t\}$
where $\Sigma_t$ is a least area plane with $\PI \Sigma_t = \Gamma_t$ by using the Lemma 3.1, i.e.
$X = \bigcup_{t=-1}^1 \Sigma_t$. Now, $S$ is a minimal surface with $\PI S=\gamma_0$. Since $S$ has
some genus, it cannot be a leave in the foliation $\{\Sigma_t\}$. So, there must be a tangential
intersection with a leaf $\Sigma_{t_1}$. This is a contradiction by maximum principle for minimal
surfaces \cite{HS}.

On the other hand, in section $4$, we gave an application of this generic uniqueness result by
applying it to the Disjoint Planes Conjecture. As described in Remark 4.1, this conjecture has many
important applications in $3$-manifold topology by combining it with the techniques developed in
\cite{Co1}. As Theorem 4.1 shows, to finish this conjecture, the only case needs to be ruled out is
when $\Gamma_1, \Gamma_2$ are not crossing each other and they both bounds more than one least area
planes, i.e. $\PI \Sigma_i^\pm =\Gamma_i$. One idea to rule out this case is to "foliate" the
region between $\Gamma_1$ and $\Gamma_2$ in $\SX$, by simple closed curves $\Gamma_t$ where $1\leq
t \leq 2$, such that $\Gamma_t\cap\Gamma_s = \Gamma_1\cap\Gamma_2$ for any $1\leq t,s\leq 2$.
Hence, one get the two families of the least area planes $\{\Sigma^+_t\}$ and $\{\Sigma^-_t\}$
where $\PI \Sigma^\pm_t = \Gamma_t$ (Otherwise, if one curve in the "foliation" $\Gamma_{t_0}$
bounds a unique least area plane, then this plane will be a barrier between $\Sigma_1^\pm$ and
$\Sigma_2^\pm$ ). Then the idea is to analyze the intersections of these families of least area
planes in order to get a contradiction by using the maximum principle for minimal surfaces.

\end{document}